\documentclass[11pt]{amsart}
\huge
\usepackage{amsthm}
\usepackage{amssymb}
\usepackage{euscript}
\newtheorem{theorem}{Theorem}
\newtheorem{theo}{Theorem}
\newtheorem{lemma}{Lemma}
\newtheorem{fact}{Proposition}
\newtheorem{cor}{Corollary}

\newtheorem{ex}{Example}
\theoremstyle{remark}

\begin{document}
\title {Quasihomogeneous $(\mathbb C^*)^k \times SL_2(\mathbb C)$-varieties
containing a finite number of orbits}
\author{E. V. Sharoyko}

\begin{abstract}
Let $G:= (\mathbb C^{\star })^k\times SL_2(\mathbb C)$ act linearly
on a vector space or its projectivisation.
We obtain an effective criterion to detect whether a number of orbits
in an orbit-closure is finite or not.
\end{abstract}
\maketitle

\section{Introduction}
In this work we operate over complex numbers $\mathbb C$.
Suppose $G:= (\mathbb C^{\star })^k\times SL_2(\mathbb C)$, $V$ is a rational
finite-dimensional $G$-module,  $\mathbb P(V)$ is its projectivisation,
$\EuScript X((\mathbb C^{\star})^k)\cong \mathbb Z^k$ is the lattice of
characters of the torus $(\mathbb C^{\star})^k$. As each $G$-module is
completely reducible, then $V=\bigoplus\limits_{i=1}^sV_i$, where $V_i$ is a
simple rational $G$-module. Let us denote by $V_{\chi_i, n_i}$ the space of
binary forms of degree $n_i$ where $G$ acts as follows:
$$
(t, \left (\begin{array}{cc}\kappa &\lambda \\\mu &\nu\end{array}\right ))
(f(x, y))=\chi_i(t)f(\kappa x+\mu y, \lambda x+\nu y).
$$
Here $n_i\in \mathbb Z_+$, $\chi_i\in \EuScript X((\mathbb C^{\star})^k)$.
Then every $V_i$ is isomorphic to some $V_{\chi_i, n_i}$.

In this paper we consider the actions $G:V$ and $G:\mathbb P(V)$. Our aim is
to obtain a criterion for detection whether a number of
orbits in an orbit-closure is finite or not. In \cite{Pa} this problem was
solved for actions $SL_2(\mathbb C):\mathbb P(V)$. In affine case there is
always a finite number of orbits as it is shown in \cite{Po}.

A $multiplicity$ is a map $e_i:\mathbb C\cup\{\infty \}\to \mathbb Z_+$
such that $e_i(\infty ):=\sum\limits_{a\in \mathbb C}e_i(a)$, and
$e_i(\infty)\le n_i$. For any $v_i\in V_i$ we have
$v_i=c_ix^{n_i-e_i(\infty)}
\prod\limits_{a\in \mathbb C}{(ax+y)^{e_i(a)}}$, and for any $v\in V$
one has $v=(c_iv_i)_{i=1}^s=(c_ix^{n_i-e_i(\infty)}
\prod\limits _{a\in \mathbb C}{(ax+y)^{e_i(a)}})_{i=1}^s$, $c_i\in \mathbb C$.
One can assume that $\forall i$ $c_i\ne 0$. Since the map
$V\to V$, $(c_iv_i)_{i=1}^s\mapsto  (v_i)_{i=1}^s$ is a $G$-automorphism
we may assume $c_i=1$ for all $i$.  Denote by $(u_i)_{i\in I}$, $I\subset
\{1, \ldots, s\}$ the vector $u\in V$  such that its $i$-th component is
equal to $u_i$ iff $i\in I$ and to 0 in other case.

Here we recall some definitions and formulate the main results of the paper.

Suppose $\langle\cdot, \cdot\rangle$ is the standard scalar product on
$\mathbb Q^{k+1}$. We take in $\mathbb Z^k\times \mathbb Z_+\subset
\mathbb Q^{k+1}$ the set of points
$(\chi_i, n_i)=\chi'_i\in\EuScript X((\mathbb C^{\star})^k\times~\mathbb C^{\star})$,
$i=1, \ldots s$.
This set is called {\it characteristic} and its points are called
{\it characteristic points}.
Let us consider the rays from zero to the characteristic points and
denote by $M$ the convex hull of these rays. One can see that $M$ is a cone.
To each face $F$ of the cone $M$ we assign the set of all characteristic points
containing in this face. This set {\it forms a face $F$}.
Consider the family of vectors $R=(r_1, \ldots r_k, p)$ with $p<0$ which
are orthogonal to a face $F$ and non-orthogonal to any face containing $F$.
If there exists a vector such that for any characteristic point $\chi'_i$ we have
$\langle\chi'_i, R\rangle\ge 0$, then the face $F$ is called {\it admissible}.
In this case we shall denote this face by $I(R)$.

\begin{theo}
The number of $G$-orbits in $\overline{Gv}$, $v\in V$ is finite iff
for each admissible face $I(R)$ of maximal dimension of the cone $M$
and for any integer-valued vector
$\beta=(\beta_i)_{i\in I(R)}$ such that
$\sum\limits_{i\in I(R)}\beta_i\chi_i=0$ the following conditions hold:
\begin{equation}
\label{afkritk}
\sum\limits_{i\in I(R)}e_i(a)\beta_i=0  \quad\forall a\in\mathbb C.
\end{equation}
\end{theo}

Proceed to the projective case. As in affine version we take in
$\mathbb Z^k\times \mathbb Z_+\subset\mathbb Q^{k+1}$  the points
$(\chi_i, n_i)=\chi'_i\in\EuScript X((\mathbb C^{\star})^k\times~\mathbb C^{\star})$
for $i=1, \ldots s$ and denote their convex by $C$.
To each face $F$ of $C$ we assign the set of all characteristic points
containing in this face. This set {\it forms an affine face $F$}.
Consider the family of vectors $R=(r_1, \ldots r_k, p)$ with $p<0$ which
are orthogonal to a face $F$ and non-orthogonal to any face containing $F$.
If there exists a vector $R=(r_1, \ldots r_k, p)$ with $p<0$ which
is orthogonal to a face $F$, non-orthogonal to any face containing $F$,
and directed into the polyhedron $C$, then the affine face $F$ is called
{\it admissible}. In this case we shall denote this face by $J(R)$.

\begin{theo}
The number of $G$-orbits in $\overline{G\langle v\rangle}$,
$\langle v\rangle\in \mathbb P(V)$ is finite iff
for each admissible face $J(R)$ of maximal dimension of the polyhedron $C$
and for any integer-valued vector $\beta=(\beta_i)_{i\in J(R)}$ such that
$\sum\limits_{i\in J(R)}\beta_i\chi_i=0$  and
$\sum\limits_{i\in J(R)}\beta_i=0$ the following conditions  hold:
\begin{equation}
\label{prkritk}
\sum\limits_{i\in J(R)}e_i(a)\beta_i=0 \quad\forall a\in\mathbb C.
\end{equation}
\end{theo}

\begin{ex}
Suppose $k=1$,
$V=V_{1,1}\oplus V_{1, 4}\oplus V_{2,4}\oplus V_{3,3}\oplus V_{4, 2}$,
$v=(v_1, x^2y^2, x^2(y-x)^2, xy(y-x), x^2)$, where $v_1$ is any linear form.
Consider the orbit-closure
$\overline{\mathbb C^*\times SL_2(\mathbb C)\langle v\rangle}\subset \mathbb P(V)$.
The faces \{B, C\} and \{C, D, E\} are admissible faces of maximal dimension.
For the face \{B, C\} all conditions of Theorem 2 hold.
For the face \{C, D, E\} the vector $\beta=(0,0,1,-2,1)$ doesn't satisfy the
condition (2) for $a=0$. Therefore, $\overline{G\langle v\rangle}$  contains
infinitely many orbits. Actually, for the curves  of the form
$\gamma(t)=(t^{-1}, \left (\begin{array}{cc}t^{-1}&-dt^{-1}\\0&t\end{array}\right ))$
with $e_i(d)=0$ we have
$\lim\limits_{t\to 0}\gamma(t)\langle v\rangle=\langle (0, 0, (d+1)^2x^4, d(d+1)x^3, x^2)\rangle$.
Let us prove that there is an infinite number of different $G$-orbits among
orbits of this form. Suppose
 $\Phi\langle v(d_1)\rangle=(t, \left (\begin{array}{cc}\kappa &\lambda \\\mu &\nu\end{array}\right ))
\langle v(d_1)\rangle=\langle v(d_2)\rangle $. Since $\langle x\rangle$
is $\Phi$-invariant then  we have $\mu=0, \nu=\kappa^{-1}$. As $\Phi x^2=x^2$
then we obtain $\kappa=\pm 1$. The following conditions are:
$t^2(d_1+1)^2=(d_2+1)^2$ ¨ $t^3d_1(d_1+1)=d_2(d_2+1)$.
It's easy to see that for any $d_1$ there exists only a finite number of $d_2$
such that $\langle v(d_2)\rangle$ is in the same orbit. Therefore, the number
of orbits in $\overline{G\langle v\rangle}$ is infinite.
\end{ex}

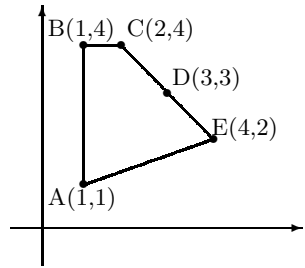
\unitlength 1mm 
\linethickness{0.4pt}
\ifx\plotpoint\undefined\newsavebox{\plotpoint}\fi 
\qquad\qquad\qquad\qquad
\begin{picture}(38.75,34.5)(0,0)
\thicklines
\put(9.75,10.75){\line(0,1){18.5}}
\put(9.75,10.75){\circle*{1}}
\put(9.75,29.25){\line(1,0){5}}
\put(9.75,29.25){\circle*{1}}
\multiput(14.75,29.25)(.033653846,-.034340659){364}{\line(0,-1){.034340659}}
\put(14.75,29.25){\circle*{1}}
\multiput(27,16.75)(-.09550562,-.03370787){178}{\line(-1,0){.09550562}}
\put(27,16.75){\circle*{1}}
\put(20.82,23){\circle*{1}}
\put(15.5,30.5){\footnotesize{C(2,4)}}
\put(21.5,24){\footnotesize{D(3,3)}}
\put(26.75,17.25){\footnotesize{E(4,2)}}
\put(5,8){\footnotesize{A(1,1)}}
\thinlines
\put(4.25,34.5){\vector(0,1){.07}}\put(4.25,0){\line(0,1){34.5}}
\put(38.75,5){\vector(1,0){.07}}\put(0,5){\line(1,0){38.75}}
\put(5,30.5){\footnotesize{B(1,4)}}
\end{picture}

Also we describe all $G$-moduli for which the orbit-closure of any orbit
in $V$ (Corollary 1)  and in  $\mathbb P(V)$ (Corollary 2) contains a finite
number of $G$-orbits.

The author is grateful to I. V. Arzhantsev for problem statement and useful
discussions.

\section{Some Lemmas}

Let $B$ be the direct product of $(\mathbb C^*)^k$ and the subgroup
of upper-triangular matrices in $SL_2(\mathbb C)$. Then $B$ is the Borel subgroup of $G$.

We shall consider $B$-orbits in $\overline{Bv}$. It is known that
$\overline{Bv}$ intersects any $G$-orbit in $\overline{Gv}$
(\cite[III.2.5, Cor. 1]{Kr}). Hence, if there is a finite number of $B$-orbits
in $\overline{Bv}$, then  there is a finite number of $G$-orbits
in $\overline{Gv}$. In the sequel we shall see that the inverse statement is
also true. Let us say that some condition holds {\it for almost all orbits}
if it holds for all orbits except finitely many.

Two following lemmas generalize Lemma and Proposition 1 from \cite{Pa}.

Let $\gamma :\mathbb C^{\star} \to B$ be a curve in $B$. We shall prove that
any $B$-orbit in $\overline{Bv}$ contains a point $\lim\limits_{t\to 0}\gamma (t)v$
for some special curve $\gamma(t)$.

\begin{lemma}
Suppose $v\in V$, $w\in \overline{Bv}$. Then there exist
$p,q,r_1, \ldots r_k\in \mathbb Z$ with $q<-p$, $c\in \mathbb C$,
and a polynomial $h\in \mathbb C[t]$ with $h(0)=-1$ and $\deg h<-p-q$ such that
$$\lim\limits_{t\to 0}(t^{r_1}, \ldots t^{r_k}, \left (\begin{array}{cc}t^p&ch(t)t^q\\0&t^{-p}\end{array}\right ))v\in Bw. $$
\end{lemma}

\begin{proof}
Let $\gamma(t)\in B(\mathbb C((t)))$ be an analitic curve in $B$ such that
$\lim\limits_{t\to 0}\gamma (t)v=w'\in Bw$.
If $\delta (t)\in B(\mathbb C((t)))$ such that
$\lim\limits_{t\to 0}\delta(t)\in B$, then
$\lim\limits_{t\to 0}(\delta (t)\cdot \gamma (t)v)=\delta (0)w'$,
and without loss of generality we can  replace $\gamma (t)$ by
$\delta (t)\gamma (t)$.

Consider $p, q, r_1, \ldots r_k\in \mathbb Z, c\in \mathbb C, f_1, \ldots f_k,
g_1, g_2\in \mathbb C[[t]],
f_1(0)\ne ~0, \ldots f_k(0)\ne ~0, g_1(0)\ne ~0, \linebreak g_2(0)=-g_1(0)$ and
$$\gamma(t)=(t^{r_1}f_1, \ldots t^{r_k}f_k,
\left (\begin{array}{cc}t^pg_1&cg_2t^q\\0&t^{-p}g_1^{-1}\end{array}\right )).$$
\newline If $á=0$ or $q\ge -p$ we put
$$\delta(t):=(f_1^{-1}, \ldots f_k^{-1},\left (\begin{array}{cc}g_1^{-1}&-t^{p+q}cg_2\\0&g_1\end{array}\right )).$$
\newline Then $\delta (t)\gamma (t)=(t^{r_1}, \ldots t^{r_k},\left (\begin{array}{cc}t^p&0\\0&t^{-p}\end{array}\right ))$.
\newline If $á\ne 0$ and $q<-p$ we put
$$\delta(t):=(f_1^{-1}, \ldots f_k^{-1}, \left (\begin{array}{cc}g_1^{-1}&ch'(t)\\0&g_1\end{array}\right )),$$
\newline where $h'\in\mathbb C[[t]]$ such that $h:=h't^{-p-q}g_1^{-1}+g_2g_1^{-1}\in \mathbb C[t]$
is a polynomial of degree less than $-p-q$. Then
$h(0)=g_2(0)g_1(0)^{-1}=-1$ and
$\newline\delta (t)\gamma (t)=(t^{r_1}, \ldots t^{r_k},\left (\begin{array}{cc}t^p&cht^q\\0&t^{-p}\end{array}\right ))$.
\end{proof}

Suppose $q<-p, \deg h<-p-q$. We shall compute $\lim\limits_{t\to 0}\gamma (t)v$
and prove that in almost all $B$-orbits in $\overline {Bv}$ a vector of
some standard form is contained.

Now we give some notations:
\newline $p_i(d)=\prod\limits_{a\in \mathbb C, a\ne d}(a-d)^{e_i(a)}$
as $d\in \mathbb C$ and $p_i(\infty)=1$;
\newline $I(d, R, A):=\{ i\mid  \langle\chi'_i, R\rangle +Ae_i(d)=0\}$,
$d\in\mathbb C\cup\{\infty\}$, $A\in \mathbb Z$
\newline (In particular, $I(R):=I(d, R, 0)$);
\newline $v(d, R, A):= (p_i(d) x^{n_i})_{i\in I(d,R, A)}$,
$v(d, R):=v(d, R, 0)$.
Vectors of the form $v(d,R)$ are called {\it standard}.

\begin{lemma}
Consider a vector of the form $v(d, R)$, where $R=(r_1, \ldots, r_k, p)$, $p<0$.
If $\forall i\in I(R)$ $e_i(d)=0$, then $v(d, R)\in \overline {Bv}$.
Almost all $B$-orbits in $\overline {Bv}$ contain a vector of this form.
\end{lemma}

\begin{proof}
We shall consider different types of the curve $\gamma$, where $\gamma$ is
of the form described in Lemma 1.
For an irreducible $G$-module $V=V_{\chi, n}$   we have
$$\gamma (t)v=t^{\langle\chi', R\rangle-pe(\infty)}x^{n-e(\infty)}
\prod\limits _{a\in \mathbb C}{((at^p+ch(t)t^q)x+t^{-p}y)^{e(a)}},$$
where $R=(r_1, \ldots, r_k, p)$, $\chi'=(\chi, n)$.

1. $p=0, c=0$.
\newline For an irreducible $G$-module $V$ we have
$$\lim\limits_{t\to 0}\gamma (t)v=
\lim\limits_{t\to 0}t^{\langle\chi', R\rangle}x^{n-e(\infty)}
\newline\prod\limits_{a\in \mathbb C}(ax+y)^{e(a)}.$$
The limit exists iff $\langle\chi', R\rangle\ge 0$.
\newline For a reducible $G$-module $V$ the limit exists iff $\forall i$
$\langle\chi'_i, R\rangle\ge 0$. In this case it is equal to $(v_i)_{i\in I(R)}$.
\newline Since $\lim\limits_{t\to 0}\gamma (t)v_i$ is either $0$ or $v_i$,
there is a finite number of the vectors of this form.

2. $p>0, c=0$.
\newline For an irreducible $G$-module $V$ we have
$$\lim\limits_{t\to 0}\gamma (t)v=\lim\limits_{t\to 0}t^{\langle\chi', R\rangle-2pe(\infty )}
\newline x^{n-e(\infty )}y^{e(\infty )}.$$
\newline The limit exists iff $\langle\chi', R\rangle-2pe(\infty)\ge ~0$.
\newline For a reducible $G$-module $V$ the limit exists iff  $\forall i$
$\langle\chi'_i, R\rangle-2pe_i(\infty )\ge 0$. In this case it is equal to
$(x^{n_i-e_i(\infty )}y^{e_i(\infty )})_{i\in I(\infty , R, -2p)}.$
\newline Since $\lim\limits_{t\to 0}\gamma (t)v_i$ is either $0$ or
$x^{n_i-e_i(\infty )}y^{e_i(\infty )}$, there is a finite number of the
vectors of this form.

3. $p<0, c=0$.
\newline For an irreducible $G$-module $V$ we have
$$\lim\limits_{t\to 0}\gamma (t)v=
\lim\limits _{t\to 0}t^{\langle\chi', R\rangle-2pe(0)}
\newline x^{n-e(0)}y^{e(0)}
\prod\limits_{a\in \mathbb C, a\ne 0}a^{e(a)}.$$
\newline The limit exists iff $\langle\chi', R\rangle-2pe(0)\ge ~0$.
\newline For a reducible $G$-module $V$ the limit exists iff $\forall i$
$\langle\chi'_i, R\rangle-2pe_i(0)\ge 0$. In this case it is equal to
$(x^{n_i-e_i(0)}y^{e_i(0)}p_i(0))_{i\in I(0 , R, -2p)}$.
\newline Since $\lim\limits_{t\to 0}\gamma (t)v_i$ is either $0$ or
$x^{n_i-e_i(0)}y^{e_i(0)}p_i(0)$, there is a finite number of the
vectors of this form.

4. $p=q, h\equiv -1$ (then we have $p<0$).
\newline For an irreducible $G$-module $V$ we have
$$\lim\limits_{t\to 0}\gamma (t)v=
\lim\limits _{t\to 0}t^{\langle\chi', R\rangle-2pe(c)}
\newline x^{n-e(c)}y^{e(c)}
\prod\limits_{a\in \mathbb C, a\ne c}(a-c)^{e(a)}.$$
\newline The limit exists iff $\langle\chi', R\rangle-2pe(c)\ge ~0$.
\newline For a reducible $G$-module $V$ the limit exists iff $\forall i$
$\langle\chi'_i, R\rangle-2pe_i(c)\ge 0$. In this case it is equal to
$(x^{n_i-e_i(c)}y^{e_i(c)}p_i(c))_{i\in I(c, R,  -2p)}$.
\newline  Since $\lim\limits_{t\to 0}\gamma (t)v_i$ is either $0$ or
$x^{n_i-e_i(c)}y^{e_i(c)}p_i(c)$, there is a finite number of the
vectors of this form with $e_i(c)\ne 0$. If for any $i$ we have $e_i(c)=0$,
then this is a standard vector $v(c, R)$.

5. $p=q, c\ne 0, h\not \equiv -1$ (then we have $p<0$).
\newline Suppose $h(t)=-1+h_lt^l+\ldots h_mt^m$.
\newline For an irreducible $G$-module $V$ we have
$$\lim\limits_{t\to 0}\gamma (t)v=
\lim\limits _{t\to 0}t^{\langle\chi', R\rangle+le(c)}x^n
 (ch_l)^{e(c)}\prod\limits_{a\in \mathbb C, a\ne c}(a-c)^{e(a)}.$$
\newline The limit exists iff $\langle\chi', R\rangle+le(c)\ge 0$.
\newline For a reducible $G$-module $V$ the limit exists iff $\forall i$
$\langle\chi'_i, R\rangle+le_i(c)\ge 0$. In this case it is equal to
$((ch_l)^{e_i(c)}p_i(c)x^{n_i})_{i\in I(c, R, l)}$.
\newline Let us act on  the vector by the element
\newline $((ch_l)^{r_1/l}, \ldots (ch_l)^{r_k/l}, \left (\begin{array}{cc}(ch_l)^{p/l}&0\\0&(ch_l)^{-p/l}\end{array}\right ))$
and
$$ ((ch_l)^{\langle\chi'_i, R\rangle/l+e_i(c)}p_i(c)x^{n_i})_{i\in I(c,R, l)}=
(p_i(c)x^{n_i})_{i\in I(c, R, l)}=v(c,R, l).$$
There is a finite number of vectors of such form with $e_i(c)\ne 0$.
If for any $i$ we have $e_i(c)=0$, then the vector is of the form $v(c, R)$.

6. $p>q, c\ne 0$.
\newline For an irreducible $G$-module $V$ we have
$$\lim\limits_{t\to 0}\gamma (t)v=
\lim\limits _{t\to 0}t^{\langle\chi', R\rangle+(q-p)e(\infty )}x^n(-c)^{e(\infty )}.$$
\newline The limit exists iff $\langle\chi', R\rangle+(q-p)e(\infty)\ge~ 0$.
\newline For a reducible $G$-module $V$ the limit exists iff $\forall i$
$\langle\chi'_i, R\rangle+(q-p)e_i(\infty )\ge 0$. In this case it is equal to
$((-c)^{e_i(\infty )}x^{n_i})_{i\in I(\infty ,R, q-p)}$.
Let us act on  the vector by the element
$$((-c)^{r_1/(q-p)}, \ldots (-c)^{r_k/(q-p)}, \left (\begin{array}{cc}(-c)^{p/(q-p)}&0\\0&(-c)^{-p/(q-p)}\end{array}\right ))$$
and  $(p_i(\infty )x^{n_i})_{i\in I(\infty, R, q-p)}=
v(\infty , R, q-p)$.

7. $p<q, c\ne 0$ (then we have $p<0$).
\newline For an irreducible $G$-module $V$ we have
 $$\lim\limits_{t\to 0}\gamma (t)v=
\lim\limits _{t\to 0}(-c)^{e(0)}t^{\langle\chi', R\rangle+(q-p)e(0)}
 x^n\prod\limits_{a\in \mathbb C, a\ne 0}a^{e(a)}.$$
\newline The limit exists iff $\langle\chi', R\rangle+(q-p)e(0)\ge ~0$.
\newline For a reducible $G$-module $V$ the limit exists iff $\forall i$
$\langle\chi'_i, R\rangle+(q-p)e_i(0)\ge 0$. In this case it is equal to
$(p_i(0)(-c)^{e_i(0)}x^{n_i})_{i\in I(0, R,q-p)}$.
Let us act on  the vector by the element
$$((-c)^{r_1/(q-p)}, \ldots (-c)^{r_k/(q-p)}, \left (\begin{array}{cc}(-c)^{p/(q-p)}&0\\0&(-c)^{-p/(q-p)}\end{array}\right ))$$
and   $(p_i(0)x^{n_i})_{i\in I(0,R,q-p)}=v(0,R,q-p)$.
\end{proof}

Let us remark that one $B$-orbit may contain more than one standard vector.

Consider the curves of the form 6 and 7 from Lemma 2. Each non-zero component
of the limit vector is  $p_i(\infty )x^{n_i}$ in case 6 and $p_i(0)x^{n_i}$
in case 7. However, each component has only a finite number of possible
values.
This means that the curves of the forms 6 or 7 have only a finite
number  of vectors $v(c, R, A)$ as limits.

Consequently, $\overline{Bv}$ contains infinitely many $B$-orbits iff
curves of the forms 4 and 5 from Lemma 2 with $e_i(c)=0$ have
$\lim\limits_{t\to 0}\gamma(t)v$ in an infinite number of different $B$-orbits.

\begin{lemma}
The number of $G$-orbits in $\overline{Gv}$ is finite iff the
number of $B$-orbits in $\overline{Bv}$ is finite.
\end{lemma}

\begin{proof}
It is known that $\overline{Bv}$ intersects each $G$-orbit in $\overline{Gv}$
(\cite [Korollar III.2.5.1]{Kr}). Hence, if the number of $B$-orbits is finite
then the number of $G$-orbits is finite. On the other hand, suppose that
the number of $B$-orbits is infinite. Then there exist infinitely many
standard vectors in different $B$-orbits. Let the element $g\in G$ move
one standard vector to another. Then $g$ must preserve $\langle x \rangle$ and
$g\in B$. However, every $G$-orbit contains at most one $B$-orbit of the
vector $v(d, R)$ and the number of $G$-orbits in $\overline{Gv}$ is also
infinite.
\end{proof}

Our aim is to check the finiteness of the number of $B$-orbits in
$\overline{Bv}$. It is sufficient to check the finiteness of the number of
$B$-orbits containing a vector of the form $v(c,R)$ such that $I(R)$ is an
admissible face. We shall now find out when  a vector $v(d_2, R_2)$ can be
obtained from $v(d_1, R_1)$ by $B$-action. Let us note that in this case one
has $I(R_1)=I(R_2)$ and thus we can consider each admissible face separately.
The number of orbits in a given orbit-closure is finite iff for each
admissible face there exists only a finite number of standard vector orbits.

Suppose that an algebraic group $H$ acts on an irreducible variety $X$.
Denote the minimal codimension of $H$-orbit in $X$ by $d(X, H)$. For an
arbitrary variety $Y$ the {\it modality} of the action $H:Y$ is
$$mod (Y, H)=\max\limits_{{X\subset Y,\atop X\ is\ irreducible\ and\
 H-invariant}}d(X, H).$$

\begin{fact}
Under the action $G:V$ for any $v\in V$ we have $mod(\overline{Gv}, G)=
mod(\overline{Bv}, B)\le 1.$
\end{fact}

\begin{proof}
Consider the space $\{(\alpha_ix^{n_i})_{i\in I(R)}\}$ for each admissible face
$I(R)$. The torus $(\mathbb C^*)^{k+1}$ acts on this space as follows:
$(\alpha_i)_{i\in I(R)}\mapsto (\chi'_i(t)\alpha_i)_{i\in I(R)}.$
In this space there is the curve $v(c, R)=
(\prod\limits_{a\in \mathbb C, a\ne c}(a-c)^{e_i(a)}x^{n_i})_{i\in I(R)}, c\in \mathbb C.$
If this curve intersects an infinite number of $(\mathbb C^*)^{k+1}$-orbits
then $mod(\overline{Bv}, B)=1$. If for all admissible faces $v(c, R)$
intesects only a finite number of orbits of the torus, then the modality is zero.

It is easy to prove that $mod(\overline{Gv}, G)=mod(\overline{Bv}, B)$
as in Lemma 3.
\end{proof}

Let us note that the inequality from Proposition 1 follows also from the paper
of E. B. Vinberg \cite{V}.

\section{Affine case}

We shall obtain a criterion to detect whether the number of $B$ orbits
in $\overline{Bv}$ is finite or not. Almost all $B$-orbits in $\overline{Bv}$
contain  $\lim\limits_{t\to 0}\gamma (t)v$ for $\gamma(t)$ of the form
described in Lemma 1 with $p=q<0, e_i(c)=0$ $\forall i$. By Lemma 2,
for almost all orbits in $\overline{Bv}$ this limit is $v(c, R)$.
We need to get the conditions which hold as we have $b\cdot v(d_1, R)=v(d_2, R)$
for some $b\in B$, $d_1, d_2\in \mathbb C$, $R=(r_1, \ldots, r_k, p)$, where
$I(R)$ is an admissible face.

We shall use the following well-known fact (see, for example, \cite{St}):

\begin{fact}
Suppose $V=\mathbb C^m$, the k-dimensional torus $T$ acts on $V$ multiplying
$i$th coordinate by the character $\chi_i$. Consider two points
$x=(x_1, \ldots x_m)$ and $y=(y_1,\ldots y_m)$ of $V$ such that
$x_1\ldots x_my_1\ldots y_m\ne 0$. Then $x$ and $y$ are in the same $T$-orbit
iff for any $\beta=~(\beta_1, \ldots \beta_m)\in \mathbb Z^m$ such that
$\beta_1\chi_1+ \ldots\beta_m\chi_m=0$ the condition
$x_1^{\beta_1}\ldots x_m^{\beta_m}=y_1^{\beta_1}\ldots y_m^{\beta_m}$
is fulfilled.
\end{fact}

\begin{theorem}
The number of $G$-orbits in $\overline{Gv}$, $v\in V$ is finite iff
for each admissible face $I(R)$ of maximal dimension of the cone $M$
and for any integer-valued vector
$\beta=(\beta_i)_{i\in I(R)}$ such that
$\sum\limits_{i\in I(R)}\beta_i\chi_i=0$ the following conditions  hold:
$$
 \sum\limits_{i\in I(R)}e_i(a)\beta_i=0 \quad \forall a\in\mathbb C.\leqno(1)
$$
\end{theorem}

\begin{proof}
By Lemma 3, it is sufficient to consider the number of $B$-orbits in
$\overline{Bv}$. We consider $b\cdot v(d_1, R)=v(d_2,R)$ and we shall obtain
the conditions which hold as there is an infinite number of $d_2$ for
some $d_1$. We have
\newline $(t, \left (\begin{array}{cc}\kappa&\lambda\\0&\kappa^{-1}
\end{array}\right ))p_i(d_1)x^{n_i}=p_i(d_2)x^{n_i}$.
\newline $\chi'_i(t, \kappa)\prod\limits_{a\in C, a\ne d_1}(a-d_1)^{e_i(a)}=
\prod\limits_{a\in C, a\ne d_2}(a-d_2)^{e_i(a)}.$
\newline Since the equation of the face $I(R)$ is $\langle\chi'_i, R\rangle=0$,
then we have $\chi'_i(t, \kappa)=\chi_i(t')$
\newline ($\chi'_i(t, \kappa)=t_1^{\alpha_{1i}}
\ldots t_k^{\alpha_{ki}}\kappa^{n_i}=t_1^{\alpha_{1i}}
\ldots t_k^{\alpha_{ki}}\kappa^{-\alpha_{1i}r_1-\ldots -\alpha_{ki}r_k}=
(t_1\kappa^{-r_1})^{\alpha_{1i}}\ldots t_k\kappa^{-r_k})^{\alpha_{ki}}=
\chi_i(t'))$.

Therefore,
$$\chi_i(t')\prod\limits_{a\in C, a\ne d_1}(a-d_1)^{e_i(a)}=
\prod\limits_{a\in C, a\ne d_2}(a-d_2)^{e_i(a)}.$$
Hence, the action is reduced to the action of the torus
and, using Proposiion 2, we obtain that a $B$-orbit contains an infinite
number of the vectors of the form $v(d,R)$ iff
$\prod\limits_{i\in I(R)}\prod\limits_{a\in \mathbb C}(a-c)^{e_i(a)\beta_i}$
is independent of $c$ for any $\beta$ such that
$\sum\limits_{i\in I(R)}\beta_i\chi_i=0$.
And it is equivalent to the following condition:
\newline for any $a\in \mathbb C$ and any $\beta$ such that
$\sum\limits_{i\in I(R)}\beta_i\chi_i=0$   we have
$\sum\limits_{i\in I(R)}e_i(a)\beta_i=0$.

Finally, any admissible face is contained in an admissible face of
maximal dimension. Hence, if condition (\ref{afkritk}) is fulfilled for
a bigger face then it is fulfilled for its subface.
So, it is sufficient to deal only with faces of maximal dimension.
\end{proof}

Now let us describe the $G$-moduli which have only a finite number of
orbits in each orbit-closure.

Suppose that $I(R)$ is an admissible face of maximal dimension. A {\it
character matrix} for this face is the matrix with the coordinates of characters
$\chi'_i$,  $i\in I(R)$ as its columns.

\begin{cor}
Any orbit-closure in $V$ contains a finite number of $G$-orbits iff
for each admissible face  $I(R)$ of the cone $M$ the characters $\chi_i$
$(i\in I(R), n_i\ne ~0)$ are linearly independent over $\mathbb Q$. Particularly,
it is true if $M$ does not contain admissible faces.
\end{cor}
\begin{proof}
By Theorem 1, $\overline{Gv}$ contains a finite number of $G$-orbits iff for
each admissible face $I(R)$ and for any $c\in \mathbb C$ the vector
$(e_i(c))_{i\in I(R)}$ is a rational linear combination of the rows
of character matrix. The vector $(e_i(c))_{i\in I(R)}$ is called the {\it vector
of mutiplicities}. Let us note that $E_j$ (the vector with 1 for $j$th
coordinate and 0 for others) is a vector of multiplicities for some $c$ iff
$j\in I(R)$ ¨ $n_j\ne 0$. On the other hand, any vector of multiplicities
is a linear combination of $E_j$.
\end{proof}

\section{Projective case}

\begin{theorem}
The number of $G$-orbits in $\overline{G\langle v\rangle}$,
$\langle v\rangle\in \mathbb P(V)$ is finite iff
for each admissible face $J(R)$ of maximal dimension of the polyhedron $C$
and for any integer-valued vector $\beta=(\beta_i)_{i\in J(R)}$ such that
$\sum\limits_{i\in J(R)}\beta_i\chi_i=0$  and
$\sum\limits_{i\in J(R)}\beta_i=0$ the following conditions hold:
$$
\sum\limits_{i\in J(R)}e_i(a)\beta_i=0 \quad\forall a\in\mathbb C.\leqno(2)
$$
\end{theorem}

\begin{proof}
To obtain the criterion in the case of projective action $G:\mathbb P(V)$,
we shall consider the linear action of bigger group $\mathbb C^*\times G$
on $V$ with characters $\hat{\chi_i}=(1, \chi_i)\in
\EuScript X((\mathbb C^{\star})^{k+1})$. Let the cone $M$ correspond to this
action.

\begin{lemma}
The number of $G$-orbits in $\overline{G\langle v\rangle}$,
$\langle v\rangle\in \mathbb P(V)$ is finite iff
for each admissible face $I(R)$ of maximal dimension of the cone $M$
and for any integer-valued vector $\beta=(\beta_i)_{i\in I(R)}$ such that
$\sum\limits_{i\in I(R)}\beta_i\chi_i=0$  and
$\sum\limits_{i\in I(R)}\beta_i=0$ the following conditions hold:
$$
\sum\limits_{i\in I(R)}e_i(a)\beta_i=0 \quad \forall a\in\mathbb C.
$$
\end{lemma}

\begin{proof}
It follows from Theorem 1.
\end{proof}

If we return from characters $\hat{\chi_i}$  back to $\chi_i$, then
a set of the form $J(R):=\{ i\mid (\chi'_i, R)=\min\limits_j{(\chi'_j, R)}\}$
corresponds to each admissible face  $I(R)$ of the cone.
This set forms an admissible affine face of the
intersection of the cone and the plane ${x_1=1}$.
\end{proof}

In the case $G=SL_2$ the polyhedron $C$ is an arc and the only its
admissible face is its right vertex. Theorem 2 shows that
$\overline {SL_2\langle v\rangle}$
contains a finite number of $G$-orbits iff all components of maximal degree
of the vector $v$ coincide.
This fact was originally proved in \cite[Prop. ~4]{Pa}.

\begin{cor}
Any orbit-closure in $\mathbb P(V)$ contains a finite number of $G$-orbits
iff for each admissible face  $J(R)$ of the polyhedron $C$ the characters
$\hat\chi_i$ $(i\in J(R), n_i\ne ~0)$ are linearly independent over $\mathbb Q$.
\end{cor}

\begin{proof}
It easily follows from Corollary 1.
\end{proof}

\begin{ex}
Consider the action $SL_2:\mathbb P(V_{n_1})\times\ldots\times\mathbb P(V_{n_m})$.
Then any orbit-closure contains only a finite number of orbits.
To prove it, we take the action of the group $SL_2\times(\mathbb C^{\star})^m$
on $V=V_{n_1}\oplus\ldots V_{n_m}$  with the character $\chi_i$ on $V_{n_i}$.
Here $\chi_1=(1,0,\ldots 0),\ldots,\chi_m=(0,\ldots 0,1)$. Now we use
Corollary 1.
\end{ex}

Chair of Higher Algebra, Faculty of Mechanics and Mathematics, Moscow
State University, Leninskie Gory, GSP-2, Moscow, 119992, Russia

e-mail: sharojko@mccme.ru

\end{document}